\documentclass[smallextended]{svjour3}  

\usepackage{amsmath,amssymb,amsfonts}
\usepackage{tikz}

\usepackage{amsmath,amssymb}
\usepackage{graphics}

\usepackage{verbatim}
\usepackage{color}

 \newcommand{\vet}{V^\star}
\newcommand{\het}{H^\star}
\newcommand{\vetn}{V^\star_\nom}
\newcommand{\hetn}{H^\star_\nom}

\def \calA{\mathcal{A}}
\def \calB{\mathcal{B}}
\def \calC{\mathcal{C}}
\def \calD{\mathcal{D}}
\def \Var{\hbox{{\rm Var}}}
\def\nom{\textup{nom}}

\newenvironment{proofof}{\noindent {\em Proof of }}{\hfill \hspace*{1pt}
\hfill $\square$}

\newcommand{\nbi}{{p}}            
\newcommand{\esp}{\mathrm{E}}              
\newcommand{\var}{\mathrm{Var}}            

\newcommand\ds{\displaystyle}

\newcommand{\RR}{\mathbb{R}}

\newcommand\real{\ensuremath{{\mathbb R}}}

\newcommand{\norm}[1]{\left\Vert#1\right\Vert}
\newcommand{\abs}[1]{\left\vert#1\right\vert}

\newcommand{\R}{\mathbb R}

\graphicspath{{figures/}{EnCours/}}
\usepackage{url}

\begin{document}

\title{Global sensitivity analysis for the boundary control of an open channel}

\author{Alexandre Janon, Ma\"elle Nodet, Christophe Prieur and Cl\'{e}mentine Prieur}
\institute{Alexandre Janon is with Department of Mathematics, Universit\'e Paris-Sud, Orsay, France.   {\tt\small Alexandre.Janon@math.u-psud.fr}, Ma\"elle Nodet and Cl\'{e}mentine Prieur are with Universit\'e Grenoble Alpes, CNRS, and INRIA, Grenoble, France. {\tt\small maelle.nodet@inria.fr},  {\tt\small clementine.prieur@imag.fr}, Christophe Prieur is with Gipsa-lab, Grenoble, France.
    {\tt\small christophe.prieur@gipsa-lab.fr}}


\maketitle

\begin{abstract}  
The goal of this paper is to solve the global sensitivity analysis for a particular control problem.
More precisely, the boundary control problem of an open-water channel is considered, where the boundary conditions are defined by the position of a down stream overflow gate and an upper stream underflow gate. The dynamics of the water depth and of the water velocity are described by the Shallow Water equations,  taking into account the bottom and friction slopes. Since some physical parameters are unknown, a stabilizing boundary control is first computed for their nominal values, and then a sensitivity analysis is performed to measure the impact of the uncertainty in the parameters on a given {\it to-be-controlled output}. The unknown physical parameters are described by some probability distribution functions. Numerical simulations are performed to measure the first-order and total 
sensitivity indices.
\end{abstract}

\section{Introduction}
In this work we consider the boundary stabilization of an open channel. The model we consider is described by the Shallow-Water equations, which are conservation laws perturbed by non-homogeneous terms due to the effects of the bottom slope, the slope's friction, and the lateral supply. The boundary actions are defined as the position of both spillways located at the extremities of the reach. In \cite{CoronBastinAndrea:SIAM:08,SantosPrieur:IEEE:08}, the authors designed stabilizing {boundary output feedback controllers}, with an exponential convergence to the equilibrium of water level and water flow. Our interest is motivated by the following remark: in a given real open channel, many of the  involved parameters are uncertain, e.g., because of measurement uncertainties (bottom slope, friction slope, \ldots). Our aim is then to study the sensitivity of the efficiency of the control of the open channel with respect to the uncertainties in these parameters.

This problem is related to the insensitizing problem, which consists in finding a control function such that  some functional of the state is locally insensitive to the perturbations of one given parameter (usually the initial condition in the literature). 
For the semilinear heat equation, we can mention  \cite{BodartFabre:JMAA:1995} for bounded domain and \cite{Teresa:ESAIM:1997} for unbounded domain (see also \cite{TeresaZuazua:CPAA:2009} for a complete insensitization). Other works include  \cite{Teresa:CPDE:2000} for insensitizing controls of semilinear parabolic equation where some data of the heat equation are incomplete and unknown (see also \cite{bodart2004existence} for the heat equation in presence of superlinear nonlinearity) and \cite{MicuOrtegaTeresa:AML:2004} when the control and the observation regions are disjoint. When focusing on hyperbolic systems and fluid mechanics, similar works have been caried out, see in particular \cite{alabau:MCSS2014:insensitizing,dager2006insensitizing} for the wave equation,  \cite{fernandez2003insensitizing} for a simplified linear ocean model, \cite{guerrero2007controllability,gueye2013insensitizing} for recent papers on Stokes and Navier--Stokes equations. 

Previous papers on this topic address a local sensitivity problem, that is they study the derivative of the quantity of interest around a given value of one parameter. The work we propose here presents a double originality. Firstly we will address simultaneously the sensitivity to all uncertain parameters, and not just one. And secondly, we propose to investigate the global sensitivity, when the parameters vary around their mean values, following prescribed probability distributions (Gaussian or uniform). Therefore, the present approach can be seen as a global analysis of the control. This will be addressed using statistical techniques.

Sensitivity analysis aims to find the most sensitive parameters, i.e. parameters whose variations have the largest impact on the output quantity. Local sensitivity analysis essentially computes the derivative of the output with respect to the parameter, at a given value of the parameter. Global (stochastic) sensitivity analysis (see e.g. \cite{saltelli-sensitivity} for a review) assumes that the parameters can vary widely, either in a given range, or around a given value. In this framework, the parameters are assumed to follow suitable probability distributions. One way of measuring sensitivities is to compute sensitivity indices, such as Sobol indices \cite{SOB1993}, which quantify the contribution of a given parameter or set of parameters to the output variance: the larger the index value, the greater the sensitivity.

These indices are in general impossible to compute exactly, and must be estimated. Classical approaches of effective computation use Monte-Carlo type methods, see the survey \cite{helton2006survey}. We can cite the FAST method (Fourier Amplitude Sensitivity Testing, \cite{cukier1978nonlinear,TIS2012}) which uses the Fourier decomposition of the output function, the polynomial chaos expansion method \cite{sudret2008global}, and the Sobol pick-freeze scheme \cite{SOB1993,sobol2001global,janon:2013,JAN2014} which uses sampled replications of model outputs. The Monte-Carlo approach requires a large number of model runs (e.g. around one thousand for one parameter). As in general the model is complex and requires large computing time, it is beneficial to replace the full model by a metamodel, that is an approximate but fast model. In this work we used the reduced basis method \cite{grepl2007efficient,grepl2005posteriori,nguyen2005reduced,veroy2005certified,janon2013m2an}, but we can also mention kriging method \cite{sant:will:notz:2003}, interpolation kernels \cite{schaback2003mathematical}, and we refer to \cite{fang2005design} for a review on metamodeling methods. The reduced basis method consists in solving the discrete model partial differential equation in a smaller dimension space, i.e. to look for a solution in a space spanned by a reduced basis instead of a large generic basis (such as finite elements). The advantage of this method is that it can provide certified error bounds that allow us to quantify the information loss between the full model and the metamodel, and therefore to provide certified Confidence Intervals (CI) for the sensitivity indices \cite{janon2011uncertainties}.

This paper is organized as follows. Section \ref{sec:2} presents the model and states the problem. Section \ref{sec:3} presents the sensitivity analysis in our context, and describes the numerical computation of indices using Monte Carlo approach and  reduced basis metamodeling. Section \ref{numeric} presents numerical results. We conclude and give outlooks in Section \ref{sec:conc}.
The Appendix collects the proof of some intermediate results. This work is an extension of the paper of the same title presented at the 2014 CDC conference \cite{janon:hal-01065886}, where no proof is given, a simpler {to-be-controlled output} is considered and less parameters are included in the sensitivity analysis.

\section{Boundary control of an open channel and problem statement}\label{sec:2}
\subsection{Quasilinear equation}

Let us consider the classical Shallow Water equations describing the flow dynamics inside of an open-channel. For an introduction of such model and related control problems see e.g. \cite{coron-book07,BastinCoronAndrea:IFAC:08,SantosPrieur:IEEE:08}. This model describes the space and time-evolution of the water depth $H=H(x,t)$ and horizontal water velocity $V=V(x,t)$ and is written as follows, for all $(x,t)\in[0,L]\times \RR_+$,
\begin{equation}
\label{nonlinear:SW}
\partial_t \left(\begin{array}{c} H\\ V\end{array}\right)
+\partial_x\left(\begin{array}{c}  HV \\\frac{1}{2}V ^2 +gH\end{array}\right)+\left(\begin{array}{c}  0\\ g (S_f -S_b)\end{array}\right) =0\ ,
\end{equation}
where $L$ stands for the length of the pool, $g$ is the gravity constant, $S_b$ is the bottom slope 
and $S_f$ is the friction slope. Moreover it is assumed that the bottom slope $S_{b}$ doesnot depend on $x$.

\begin{figure}[htbp]
\begin{center}
\def\svgwidth{8cm}
\input{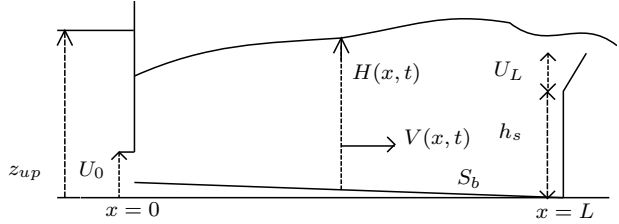}
\caption{Sketch of a channel: one reach with a downstream overflow gate and an upperstream underflow gate.\label{ft1}} \end{center}
\end{figure}

Let us denote the water flow by $Q$. It is given by $Q(x,t)=B H(x,t) V (x,t)$ where $B$ is the channel width. 
In the present work,
we suppose there are two gates, one at $ x=0$ and one at $x=L$, which 
are respectively
\begin{itemize}
\item {a submerged underflow gate}:
 \begin{eqnarray}
 \label{tau3} 
Q(0,t)&=&U_{0}B \mu_{0}\sqrt{2g(z_{up}-H(0,t))},
  \end{eqnarray}
  where  $z_{up}$ is the water level before the gate,  $\mu_0$ the water flow coefficient and $U_{0}$ the position of the spillway (see Fig. \ref{ft1}),
\item {a submerged overflow gate}:
\begin{eqnarray} \label{under}
H(L,t)&=&\left(\frac{Q^2(L,t)}{2gB ^2\mu_L^2}\right)^{1/3}+h_{s}+U_L,
\end{eqnarray}
where $h_s$ is the height of the fixed part of the overflow gate, $\mu_L$ the water flow coefficient at this gate and $U_{L}$ the position of the spillway (see Fig. \ref{ft1}).
\end{itemize}
The controls are the positions $U_0$ and $U_L$ of both spillways located at the extremities of the pool and related to the state variables $H$ and $Q$.

{Other kinds of boundary conditions (BC)
can be considered e.g., two submerged underflow gates (or two submerged
overflow gates) at $x=0$ and at $x=L$.}

There are sufficient stability conditions written in terms of the boundary conditions for the stability of (\ref{nonlinear:SW}) with the boundary conditions (\ref{tau3})-(\ref{under}). These sufficient conditions exploit Lyapunov functions (see e.g. \cite{CoronBastinAndrea:SIAM:08}), or analysis of the characteristic curves (see \cite{PrieurWinkinBastin08,libook:94}).

There are various empirical models that are available in the literature for the friction slope  (see e.g. \cite{BastinCoronAndrea:NHM:09,SantosPrieur:IEEE:08} which are considering two different models). 
%
Let us pick the model of \cite{BastinCoronAndrea:NHM:09} for $S_f$, that is 
\begin{equation}
\label{friction:slope}
S_f =
C \frac{V^2}{H},
\end{equation}
where $C$ is a constant friction coefficient.

\subsection{Linearized model and stabilizing control laws}
A steady-state solution or \emph{equilibrium} of (\ref{nonlinear:SW}) is a time-independent solution of this equation. Let us consider a space-independent equilibrium, denoted $H^{\star}, V^{\star}$ and defined by:
\begin{equation}\label{e:defvtoilhtoil} V^ \star= \left( \frac{S_b Q^\star}{BC} \right)^{1/3}, \;\; H^\star = \frac{Q^\star}{BV^\star} . \end{equation}
Using (\ref{nonlinear:SW}) and (\ref{friction:slope}), it implies $S_f = S_b$ and $S_b H^\star =C (V^\star) ^2$. The linearized Shallow Water equations around such an equilibrium are computed in \cite{BastinCoronAndrea:NHM:09}. Denoting the deviations of the state with respect to such an equilibrium by $h(x,t)=H(x,t)-H^\star$ and $v(x,t)=V(x,t)-V^\star$, we have:
$$
\partial_t \left(\begin{array}{c} h \\v \end{array}\right)+\left(\begin{array}{cc} V^\star& H^\star\\ g & V^\star\end{array}\right)\partial_x\left(\begin{array}{c} h \\v \end{array}\right)
+ \left(\begin{array}{cc} 0 & 0 \\ - \frac{g S_b}{H^\star} & \frac{2 gS_b}{V^\star}\end{array}\right)\left(\begin{array}{c} h \\v \end{array}\right)= 0  \ ,
$$
or equivalently, by recalling that $S_b H^\star =C (V^\star) ^2$,
\begin{equation}
\label{linear:SW}
\partial_t \left(\begin{array}{c} h \\v \end{array}\right)+\left(\begin{array}{cc} V^\star& H^\star\\ g & V^\star\end{array}\right)\partial_x\left(\begin{array}{c} h \\v \end{array}\right)
+ \left(\begin{array}{cc} 0 & 0 \\ - \frac{g C (V^\star)^2}{(H^\star)^2} & \frac{2 g C V^\star}{H^\star}\end{array}\right)\left(\begin{array}{c} h \\v \end{array}\right)= 0  \ .
\end{equation}
We now introduce the classical characteristic coordinates defined by:
\begin{equation}
\label{charac:coord}
\begin{array}{c}
\xi _ 1(x,t) =v(x,t) + h(x,t) \sqrt{\frac{g}{H^\star}}
\smallskip\\
\xi _ 2(x,t) =v(x,t) - h(x,t) \sqrt{\frac{g}{H^\star}}
\end{array}
\end{equation}
for all $(x,t)\in [0,L]\times \RR_{+}$, and the characteristic velocities:
\begin{equation}
\label{charac:veloc}
\lambda_ 1=V^\star + \sqrt{gH^\star}
\ , \; 
- \lambda_ 2=V^\star - \sqrt{gH^\star}\ .
\end{equation}
Assume that the flow is fluvial, that is 
$
gH^\star > {V^\star}^2
$.
Under this condition, the characteristic velocities have opposite signs, that is $\lambda_{1}>0$ and $-\lambda_{2}<0$.

The linearized Shallow Water equations (\ref{linear:SW}) may be rewritten as
\begin{equation}
\label{linear:SW:charac}
\partial_t \left(\begin{array}{c} \xi_1 \\\xi_2 \end{array}\right)+\left(\begin{array}{cc} \lambda _1& 0\\  0& -\lambda _2 \end{array}\right)\partial_x \left(\begin{array}{c} \xi_1 \\\xi_2 \end{array}\right)\\
+ \left(\begin{array}{cc} \gamma  & \delta \\ \gamma  & \delta\end{array}\right)\left(\begin{array}{c} \xi_1  \\\xi_2 \end{array}\right)= 0 
\end{equation}
with the parameters
\begin{equation}\label{eq:def:gamma:delta}
\begin{array}{c}
\gamma = gC \frac{(V^\star)^2}{H^\star}\left(
\frac{1}{V^\star}- \frac{1}{2\sqrt{gH^\star}} \right) 
\ , \smallskip\\
\delta = gC \frac{(V^\star)^2}{H^\star}\left(
\frac{1}{V^\star}+ \frac{1}{2\sqrt{gH^\star}} \right) 
\ .
\end{array}\end{equation}

We then have the following proposition (see Appendix \ref{appendix:A} for the proof which is inspired 
by \cite{SantosPrieur:IEEE:08,BastinCoronAndrea:NHM:09}):

\begin{proposition}
\label{def:controls:lemma}
Given any constant values $k_{0}$ and $k_L$, let us define the controls $U_0$ and $U_L$ by, for all $t\geq 0$,
\begin{equation}
\label{def:U0}
U_0 (t) = \frac{H(0,t)\left( V ^\star -\frac{1+ k_0}{1- k_0} ( H (0,t) - H ^\star ) \sqrt{\frac{g}{H^\star}} \right)}{\mu _0\sqrt{2g(z_{up}-H(0,t))}} ;
\end{equation}

\begin{equation}
\label{def:UL}
U_L(t) = -\left(
\ds\frac{H (L,t)\left(V ^\star  + \frac{1+ k_L}{1- k_l} ( H( L,t) - H ^\star) \sqrt{\frac{g }{H^\star}}
\right)
}{\sqrt{2 g} \mu_L }
\right)^{\frac{2}{3}} +H(L,t) - h_s.
\end{equation}
Then the boundary conditions (\ref{tau3}) and (\ref{under}) may be rewritten as
 \begin{equation}
\label{bound:linear:SW}
\left(\begin{array}{c} \xi_1 (0,t) \\\xi_2 (L,t) \end{array}\right) =
\left(\begin{array}{cc} 0& k_{0}\\k_{L}&0 \end{array}\right) \left(\begin{array}{c} \xi_1 (L,t) \\\xi_2 (0,t) \end{array}\right) \ .
\end{equation}
\end{proposition}
Note that by defining the to-be-measured outputs as the water heights at both extremities of the channel, namely $H(0,t)$, and $H(L,t)$, the previous result defines output feedback controllers $U_0$ and $U_L$. 
Moreover, note that this is an exact expression (not coming from an approximation or a linearization of the boundary conditions (\ref{tau3}) and (\ref{under})).

It is possible to combine the previous results with Lyapunov techniques to compute stabilizing controllers. This is one of the contributions of \cite{BastinCoronAndrea:NHM:09} which is recalled here:
\begin{proposition}\label{propo:stab:control} (\cite{BastinCoronAndrea:NHM:09})
For any $(k_0,k_L)\in \RR$ such that
\begin{equation}
\label{stab:cond:bastin}
\max \left\{ |k_0| \sqrt{\frac{\lambda_1 \gamma}{\lambda_2 \delta}}, |k_L| \sqrt{\frac{\lambda_2 \delta}{\lambda_1 \gamma}} \right\} < 1 \ ,
\end{equation}
where $\gamma$ and $\delta$ are defined in (\ref{eq:def:gamma:delta}), defining $U_0$ and $U_L$ with Proposition \ref{def:controls:lemma},  
the system (\ref{linear:SW:charac}) with the boundary conditions
(\ref{bound:linear:SW}) is exponentially stable (in $L^2$-norm). More precisely, there exist $\nu > 0$ and 
$M>0$ such that, for every initial condition $(\xi_1^0,\xi_2^0) \in L^2((0, L); \RR^2)$, the solution to the Cauchy problem (\ref{linear:SW:charac}) with the boundary conditions
(\ref{bound:linear:SW}) and the initial condition
\begin{equation}
\label{xi:initial:cond}
(\xi_1 (x,0),\xi_2 (x,0)) = (\xi_1 ^0(x),\xi_2^0 (x))  \ ,\; \forall x \in (0,L) 
\end{equation}
satisfies
$$
\| (\xi_1(\cdot,t),\xi_2 (\cdot,t))\|_{L^2((0,L);\RR^2)} \leq Me^{-\nu t} \| (\xi_1^0,\xi_2^0 )\|_{L^2((0,L);\RR^2)}
.
$$
\end{proposition}

The proof of the previous result is given in \cite{BastinCoronAndrea:NHM:09} by noting that (\ref{stab:cond:bastin}) is equivalent to \cite[Condition (9)]{BastinCoronAndrea:NHM:09}.

\subsection{Problem statement}
Assume now that we want to use this study to compute a stabilizing control for a real-life channel. In this case, a legitimate question to ask is whether the controls computed using the theoretical model are accurate enough to stabilize the real-life channel according to Proposition \ref{propo:stab:control}. More specifically, most of the physical parameters in this problem are measured quantities, potentially endowed with measurement errors. Thus the theoretical channel may differ from the real-life one. And as the control has been designed on the theoretical channel, it may lead to a difference in the quality of the stabilization of the real-life (or uncertain) channel. In this paper, we want to find which physical parameters have the largest influence on this quality of stabilization. 

Indeed, if there are measurement errors, the real-life model is still governed by (\ref{linear:SW:charac}), but the true values for the parameters are unknown, {\em a priori} different from the nominal values used to model the theoretical channel. Moreover some parameters in (\ref{nonlinear:SW}) are not well modeled, e.g., the friction slope $S_f$, given by (\ref{friction:slope}), for which others models exist in the literature (see e.g. \cite{BastinCoronAndrea:NHM:09} and \cite{SantosPrieur:IEEE:08}). The only way to design the control is to use for each parameter its nominal value. It leads to nonlinear boundary conditions, which we linearize in order to keep the 
resolution simple. Let us denote with $nom$ in subscript the nominal value of a parameter. For instance, $z_{up,\nom}$ is the nominal value of $z_{up}$. The quantities $V^\star_\nom$ and $H^\star_\nom$ are defined by \eqref{e:defvtoilhtoil}, with all parameters fixed to their nominal values.
For sake of conciseness, we omit the time variable: for instance $H(0)$ stands for $H(0,t)$ for each $t$. We obtain the following proposition
(see Appendix \ref{appendix:B} for the proof):
\begin{proposition}
\label{def:controlsreal_life:lemma}
Given any constant values $k_{0}$ and $k_L$, let us define the controls $U_0$ and $U_L$ by, for all $t\geq 0$,
\begin{eqnarray}\label{def:U0rl}
U_0 (t)&=&\ds\frac{H(0,t)\left( V_\nom ^\star -\frac{1+ k_0}{1- k_0} ( H (0,t) - H_\nom ^\star ) \sqrt{\frac{g}{H_\nom^\star}} \right)}{\mu _0\sqrt{2g(z_{up,\nom}-H(0,t))}}
\medskip\\
\nonumber
U_L(t)&=&-\left(
\ds\frac{H (L,t)\left(V_\nom ^\star  + \frac{1+ k_L}{1- k_L} ( H( L,t) - H _\nom^\star) \sqrt{\frac{g }{H_\nom^\star}}
\right)
}{\sqrt{2 g} \mu_L }
\right)^{\frac{2}{3}}\\
\label{def:UL:2}&& +H(L,t) - h_{s,\nom}.
\end{eqnarray}
Then the boundary conditions (\ref{tau3}) and (\ref{under}) for the real-life model(\ref{linear:SW:charac}) are linearized as:
\begin{equation}\label{bound:expmodel1} \left( 1-\frac{\calB+\sqrt{g/\het}}{2\sqrt{g/\het}} \right) \xi_1(0,t) + \frac{ \calB+\sqrt{g/\het}}{2\sqrt{g/\het}}\xi_2(0,t)=\calA \end{equation}
\begin{equation}\label{bound:expmodel2} - \frac{ \calD-\sqrt{g/\het} }{2\sqrt{g/\het}} \xi_1(L,t) + \left( 1 + \frac{ \calD-\sqrt{g/\het} }{2\sqrt{g/\het}} \right) \xi_2(L,t)=\calC \end{equation}
where $\calA$, $\calB$, $\calC$ and $\calD$ are values defined in Appendix \ref{appendix:B}.
\end{proposition}
\begin{remark}
We remark that if the nominal parameters coincide with the true ones, Proposition \ref{def:controls:lemma} and Proposition \ref{def:controlsreal_life:lemma} provide the same boundary conditions.
\end{remark}
{Let us define by ${ \mu}=\left( h_s,B,S_b,C,z_{up},\xi_1^0, \xi_2^0,\mu_0,\mu_L\right)$ the vector of uncertain parameters. The uncertainty on these parameters is modeled by random variables. The measurements of these parameters can then be considered as random variables evaluations. Parameters $h_s$, $B$, $S_b$, $z_{up}$, $\mu_0$ and $\mu_L$  are modeled by Gaussian distributions whose means are the nominal values of the parameters and whose standard deviations reflect the uncertainties on these measurements or evaluations. Initial conditions $(\xi_1^0,\xi_2^0)$ and parameter $C$ (defining the friction slope in (\ref{friction:slope})) are modeled by uniform distributions. We refer to Table \ref{tableparams} for more details on these distributions. 
\begin{table}[htbp]
\[
{\begin{array}{l|l|l}
\textrm{Name} & \textrm{Nominal value} & \textrm{Comment} \\ \hline
h_s & 4m &   \textup{uncertainty:  } \mathcal{N}\left(4,0.03 \right)\\  
B & 80m &  \textup{uncertainty } \mathcal{N}\left(80,1.03 \right)\\   
S_b & 0.0002 &  \textup{uncertainty } \mathcal{N}\left(2 \times 10^{-4}, 2.5 \times 10^{-6}\right) \\ 
C & 0.001 &   \textup{uncertainty } \mathcal{U}\left([9 \times 10^{-4}, 0.0011]\right)\\ 
z_{up} & 10 m &  \textup{uncertainty } \mathcal{N}\left(10,0.13\right) \\  
\xi_1^0 & 0 & \textrm{initial value, uncertainty } \mathcal{U}\left([-0.01, 0.01]\right)  \\  
\xi_2^0 & 0 & \textrm{initial value, uncertainty } \mathcal{U}\left([-0.01, 0.01]\right)\\
\mu_0 & 0.65 & \textup{uncertainty } \mathcal{N}\left(0.65,0.0066\right)\\  
\mu_L & 0.65 &\textup{uncertainty }\mathcal{N}\left(0.65,0.0066\right)
\\  
\hline  
k_0 & 0.6 & \textup{known}\\          
k_L & 0.7 & \textup{known}\\ 
Q^\star & 50 & \textup{known }\ \\   
g & 9.81 & \textrm{acceleration of gravity, known} \\  
\end{array}}
\]
\caption{Uncertainty on input parameters. $\mathcal{N}(m,\sigma)$ is a normal distribution of mean value $m$ and standard deviation $\sigma$, and $\mathcal{U}([a,b])$ is the uniform distribution on $[a,b]$.} \label{tableparams}
 \end{table}

The stability of the system is then measured by the so-called {\it to-be-controlled output}:
\begin{equation}\label{e:defsortie}f(\mu)=\sqrt{ \int_{t=0}^{T^{\star}}\int_{x=0}^{L} \xi_1(x,t)^2+\xi_2(x,t)^2 \, \textup{d}x\, \textup{d}t} \, , \end{equation}
where $T^{\star}$ is a given time horizon. In our study, the parameters $(k_0,k_L)$ will be fixed by the controller. Recall that $\xi$ is governed by Equations (\ref{linear:SW:charac}) with boundary conditions given by Proposition \ref{def:controlsreal_life:lemma}. The sensitivity of the {\it to-be-controlled output} to the input parameters $\mu$ will be derived by performing a global sensitivity analysis whose bases are recalled in the next section.
 The closed-loop system is sketched out
in Figure \ref{fig:nominal}, and in Figure \ref{fig:uncertain} where the uncertainties appear.
\begin{figure}[ht]
\centerline{\input{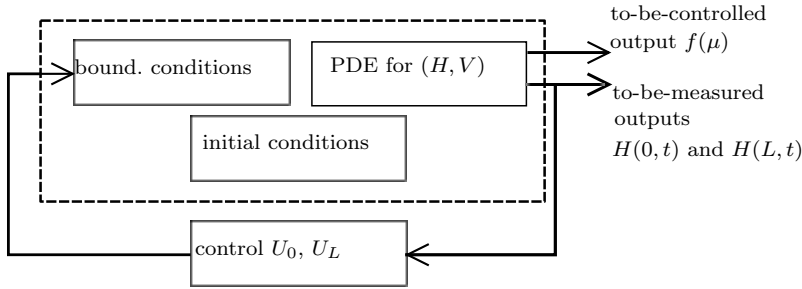}}
\caption{Control loop in case when are no uncertainties.}\label{fig:nominal}
\end{figure}

\begin{figure}[ht]
\centerline{\input{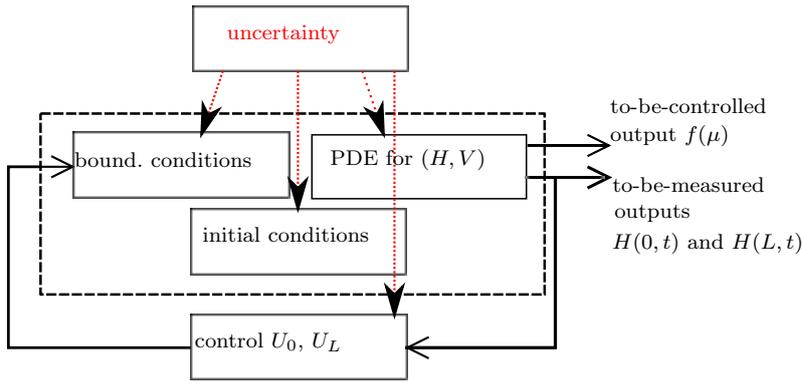}}
\caption{Control loop in case of uncertainties.}\label{fig:uncertain}
\end{figure}

\section{Sensitivity analysis}\label{sec:3}
\label{section:sa}

In the following, one wants to measure the sensitivity of the {\it to-be-controlled output} $f(\mu)$ with respect to the uncertainties on the parameter vector $\mu$. 

\subsection{Global sensitivity: a variance-based approach}
We adopt a stochastic framework. Input parameters $\mu_1$, \ldots, $\mu_p$ (here $p=9$) are assumed to be independent and are thus modeled by one-dimensional distributions, as detailed in Table \ref{tableparams}. The {\it to-be-controlled output} can then be considered as a scalar random variable $Y=f(\mu)$.

The conditional expectation $\esp(Y | \mu_i)$ is a random variable which gives the mean of $Y$ over the distributions of the $\mu_j$ ($j\not=i$), when $\mu_i$ is fixed. It is the best approximation in the mean square sense of $Y$ which depends on $\mu_i$ only. Its variance quantifies the influence of $\mu_i$ on the dispersion of $Y$. The Sobol' sensitivity indices are obtained by normalizing this variance, by the total variance of the to-be-controlled output, $Y$.
Thus, the Sobol' first order sensitivity index of input parameter $\mu_i$ is defined by
\begin{equation*}
  {\mathrm{S}}_{\{i\}} = \frac{\var\big(\esp(Y | \mu_i)\big)}{\var(Y)}.
\end{equation*}
It belongs to the interval $[0,1]$.

More generally, one can define sensitivity indices of any order $r \in \{1, \ldots , \nbi \}$, starting from the functional Analysis Of Variance (ANOVA) 
decomposition (see \cite{HOE1948} and \cite{EFS1981}). Let us first introduce some notations. We assume that $f$ is a real square integrable function, $\mathfrak{u}$ is a subset of $\{1,\ldots ,\nbi\}$, $\mathfrak{u}^c$ stands for its complement, its cardinal is denoted $r=|\mathfrak{u}|$, and  $\mathbf{\mu}_{\mathfrak{u}}$ represents the random vector with
components $\mu_i$, $i\in\mathfrak{u}$. The ANOVA decomposition
then states that 
$Y=f(\mu)$ can be uniquely decomposed into summands of increasing dimensions 
\begin{equation}
 f(\mu)=\sum_{\mathfrak{u}\subseteq \{1,\ldots ,\nbi\}} f_{\mathfrak{u}}(\mathbf{\mu}_\mathfrak{u})
\label{hoef}\end{equation}
where $f_\emptyset=\mathbb{E}[Y]$ and the other components have zero mean value and are mutually uncorrelated.
The Sobol' index  \cite{SOB1993}  of order $r=|\mathfrak{u}|$ with respect to the combination of all the variables in $\mathfrak{u}\subseteq\{1,\ldots ,p\}$ is then defined as 
\begin{equation}\label{indsob} {\mathrm{S}}_\mathfrak{u} = \frac{\sigma_\mathfrak{u}^2}{\sigma^2}=\frac{\textrm{Var}\big[f_\mathfrak{u}
(\mathbf{\mu}_\mathfrak{u})\big]}{\textrm{Var}[Y]}\, .\end{equation}

The main effect of the $i^{\textrm{th}}$ factor is measured by ${\mathrm{S}}_{\{i\}}$. Then, for $i \neq j$, the interaction effect due to the $i^{\textrm{th}}$ and 
the $j^{\textrm{th}}$ factors, that cannot be explained by the sum of the individual effects of $\mu_i$ and $\mu_j$, is measured by $ {\mathrm{S}}_{\{i,j\}}$, and 
so on (see \cite{saltelli-sensitivity}).
For any $i \in \{1,\ldots ,\nbi\}$, we also define a total sensitivity index $\mathrm{S}_{\{i\}}^{\textup{tot}}$ \cite{HOS1996} to express
the overall output sensitivity to an input $\mu_i$ by

\begin{equation}\label{sobtot}
\mathrm{S}_{\{i\}}^{\textup{tot}}=\sum_{\mathfrak{v} \subset \{1, \ldots p\} {\textrm{ such that }}i \in \mathfrak{v}} \mathrm{S}_{\mathfrak{v}} \, .
\end{equation}

For practical purposes, we also define the closed Sobol' index of order $r=|\mathfrak{u}|$ with respect to the combination of all the variables in $\mathfrak{u}\subseteq\{1,\ldots ,p\}$ as 
\begin{equation}\label{closed}
\mathrm{S}_{\mathfrak{u}}^{\textup{closed}}=\sum_{\mathfrak{v} \subset \mathfrak{u}} \mathrm{S}_{\mathfrak{v}} \, .
\end{equation} 

Let us remark that for first-order indices, ${\mathrm{S}}_{\{i\}}^{\textup{closed}}$ and ${\mathrm{S}}_{\{i\}}$ coincide.

\subsection{An estimator for Sobol' indices}\label{sobolEst}
In our context, no analytical formula is available for the Sobol' indices, which we thus need to estimate.
In this subsection, we introduce the classical Monte Carlo estimator of Sobol' indices first introduced in \cite{SOB1993}. 

We first need some notation. {Let $\mathbb{P}_{\mu}$ be the probability distribution of the random vector $\mu$.} Let $\mathfrak{u}$ be a non-empty subset of $\{1,\ldots \nbi\}$. For any $i \in \{1, \ldots, p\}$, let $\mu_i^{j,1}$ and $\mu_i^{j,2}$, $j=1, \ldots, n$ be two {\it independent and identically distributed} (i.i.d.) samples of size $n$ of the parameter $\mu_i$. Recall that the $p$ parameters $\mu_1$, \ldots, $\mu_p$ are distributed according to distributions given in Table \ref{tableparams}. We now define
$$\begin{array}{rcl}
\mathbf{\mu}_\mathfrak{u}^j & = & (\mu_i^{j,1} \, , \; i \in \mathfrak{u}) \\
\mathbf{\mu}_{\mathfrak{u}^c}^{j,1} & = &(\mu_i^{j,1} \, , \; i \in \mathfrak{u}^c) \\
\mathbf{\mu}_{\mathfrak{u}^c}^{j,2} & = & (\mu_i^{j,2} \, , \; i \in \mathfrak{u}^c) \, .
\end{array}$$
%

Finally, for $k=1$ and $2$, consider
\begin{equation}
Y_\mathfrak{u}^{j,k} = f(\mathbf{\mu}_\mathfrak{u}^j,\mathbf{\mu}_{\mathfrak{u}^c}^{j,k}).
\end{equation}

We define the estimator ${\hat{\mathrm{S}}^{\textup{closed}}}_{\mathfrak{u},n}$ of $\mathrm{S}_{\mathfrak{u}}^{\textup{closed}}$ as

\begin{equation}\label{estclas2}
{\hat{\mathrm{S}}^{\textup{closed}}}_{\mathfrak{u},n}= \frac{\displaystyle\frac{1}{n}\sum_{j=1}^n Y_\mathfrak{u}^{j,1} 
Y_\mathfrak{u}^{j,2}-\left(\frac{1}{n}\sum_{j=1}^n \frac{Y_\mathfrak{u}^{j,1} +Y_\mathfrak{u}^{j,2}}{2}\right)^2}{\displaystyle\frac{1}{n}\sum_{j=1}^n \frac{(Y_\mathfrak{u}^{j,1})^2 +(Y_\mathfrak{u}^{j,2})^2}{2}-\left(\frac{1}{n}\sum_{j=1}^n \frac{Y_\mathfrak{u}^{j,1} +Y_\mathfrak{u}^{j,2}}{2}\right)^2}\end{equation}
with $\mathfrak{u}=\{i\}$ for first-order indices and $\mathfrak{u}=\{i,j\}$ for closed second-order indices. This estimator was first introduced in \cite{Monod2006}. 
The asymptotic properties of ${\hat{\mathrm{S}}}_{\mathfrak{u},n}^{\textup{closed}}$ are stated in \cite[Propositions 2.2 and 2.5]{janon:2013}. This estimator requires a large number $n$ of model evaluations. Therefore, 
to reduce the cost, we opt for a metamodel approach, replacing equations (\ref{linear:SW:charac}) by a reduction of a linear version of (\ref{linear:SW:charac}). We refer to Section \ref{redpart} for more details on the reduction procedure.
Let us now define the respective estimators of ${\mathrm{S}}_{\{i\}}$ and ${\mathrm{S}}_{\{i\}}^{\textup{tot}}$
as
$${\hat{\mathrm{S}}}_{\{i\},n}={\hat{\mathrm{S}}}_{\{i\},n}^{\textup{closed}} \, , \;
{\hat{\mathrm{S}}}_{\{i\},n}^{\textup{tot}}=1-{\hat{\mathrm{S}}}_{\{i\}^c,n}^{\textup{closed}}\ .$$
Thanks to \cite{janon:2013} we can state the following theorem, which gives asymptotic confidence intervals for the Sobol' indices:
{\begin{theorem}
\label{ICasymp}
Assume that $\mathbb{E}(Y^4)<\infty$. Let $\alpha \in (0,1)$ (typically $\alpha=0.05$ or $0.10$). Let $\mathfrak{u} \subset \{1, \ldots , p\}$ and $i \in \{1, \ldots , p\}$.
\begin{enumerate}
\item[(i)] Let $v_{\mathfrak{u}}$ be defined by 
\begin{equation}\label{limitvar}
v_{\mathfrak{u}}^2=\displaystyle \frac{\Var \left((Y-\mathbb{E} Y)(Y_{\mathfrak{u}}-\mathbb{E}Y)-\frac{\mathrm{S}_{\mathfrak{u}}}{2}\left((Y-\mathbb{E}Y)^2+(Y_{\mathfrak{u}}-\mathbb{E}Y)^2\right)\right)}{\left(\Var Y\right)^2}\ .
\end{equation}
Then, for any consistent estimator $\hat{v}_{\mathfrak{u}}$ of $v_u$, an asymptotic confidence interval (CI) of level $1- \alpha$ for $\mathrm{S}_{\mathfrak{u}}^{\textup{closed}}$ is given by
\begin{equation*}\label{confint}
\left[ {\hat{\mathrm{S}}}_{\mathfrak{u},n}^{\textup{closed}} - z_{1- \frac{\alpha}{2}}\frac{\hat{v}_{\mathfrak{u}}}{\sqrt{n}},{\hat{\mathrm{S}}}_{\mathfrak{u},n}^{\textup{closed}} + z_{1- \frac{\alpha}{2}}\frac{\hat{v}_{\mathfrak{u}}}{\sqrt{n}}\right],
\end{equation*}
with $z_{1- \frac{\alpha}{2}}$ is the $1-\frac{\alpha}{2}$ quantile of the $\mathcal{N}(0,1)$;
\item[(ii)] Let $v_{\sim i}^2$ be defined by (\ref{limitvar}) with 
$\mathfrak{u}= \{1, \ldots i-1, i+1, \ldots , p\}$.
Then, for any consistent estimator $\hat{v}_{\sim i}$ of $v_{\sim i}$, 
an asymptotic CI of level $1- \alpha$ for $\mathrm{S}_{\{i\}}^{\textup{tot}}$ 
is given by
$$\left[ {\hat{\mathrm{S}}}_{\{i\},n}^{\textup{tot}} - z_{1- \frac{\alpha}{2}}\frac{\hat{v}_{\sim i}}{\sqrt{n}},{\hat{\mathrm{S}}}_{\{i\},n}^{\textup{tot}} + z_{1- \frac{\alpha}{2}}\frac{\hat{v}_{\sim i}}{\sqrt{n}}\right],
$$
with $z_{1- \frac{\alpha}{2}}$ is the $1-\frac{\alpha}{2}$ quantile of the $\mathcal{N}(0,1)$ distribution. \end{enumerate}
\end{theorem}}
{\begin{proofof}{{\em Theorem \ref{ICasymp}}.}
  From \cite[Proposition 2.2]{janon:2013}, we get that
$$\sqrt{n}\left(\hat{\mathrm{S}}^{\textup{closed}}_{\mathfrak{u},n}-{\mathrm{S}}^{\textup{closed}}_{\mathfrak{u}}\right) \xrightarrow[n \rightarrow + \infty]{\mathcal{D}}\mathcal{N}\left(0,v^2_{\mathfrak{u}}\right) $$
with $v_{\mathfrak{u}}^2$ defined by (\ref{limitvar}). Then,  using Slutsky's Lemma, 
we get the result in Item \textit{(i)} of Theorem \ref{ICasymp} for $\mathrm{S}_{\mathfrak{u}}$. Moreover, note that ${\mathrm{S}}_{\{i\}}^{\textup{tot}}=1-{\mathrm{S}}_{\sim i}^{\textup{closed}}$ with $\sim i = \{1, \ldots i-1, i+1, \ldots , p\}$. We thus get the result in Item \textit{(ii)} for ${\mathrm{S}}_{\{i\}}^{\textup{tot}}$.
\end{proofof}}

\medskip

{Assume now that we approximate the true model $f$ by a metamodel $\tilde{f}$. In the following we will consider reduced basis metamodeling (see Section \ref{redpart} below). 
For $k=1$ and $2$, we now consider
\begin{equation}
\widetilde{Y}_\mathfrak{u}^{j,k} =\tilde{f}(\mathbf{\mu}_\mathfrak{u}^j,\mathbf{\mu}_{\mathfrak{u}^c}^{j,k})~\label{fin}.
\end{equation}
We define $\widetilde{\mathrm{S}}^{\textup{closed}}_{\mathfrak{u},n}$ as in (\ref{estclas2}) by replacing $Y$ by $\widetilde{Y}$. Using the results in \cite{janon:2013} we have the following confidence intervals for the Sobol indices, using now the estimator based on the metamodel:
\begin{corollary}\label{metmet}
Assume that the metamodel $\tilde{f}$ depends on the Monte-Carlo sample size $n$, that there exists $c \in \mathbb{R}$ such that $\tilde{f}_n(\mu)-f(\mu) \xrightarrow[n \rightarrow + \infty]{\mathbb{L}^2(\mathbb{P}_{\mu})} c$, and that $n \, \Var \left(\tilde{f}_n(\mu)-f(\mu)\right) \xrightarrow[n \rightarrow + \infty]{}0$.
Assume that $\mathbb{E}(Y^4)<\infty$. Let $\alpha \in (0,1)$ (typically $\alpha=0.05$ or $0.10$). Let $\mathfrak{u} \subset \{1, \ldots , p\}$ and $i \in \{1, \ldots , p\}$. Then,
\begin{enumerate}
\item[(i)] for any consistent estimator $\hat{v}_{\mathfrak{u}}$ of $v_u$, where $v_{\mathfrak{u}}^2$ is defined by (\ref{limitvar}), an asymptotic CI of level $1- \alpha$ for $\mathrm{S}_{\mathfrak{u}}^{\textup{closed}}$
is given by
$$\left[ {\widetilde{\mathrm{S}}}_{\mathfrak{u},n}^{\textup{closed}} - z_{1- \frac{\alpha}{2}}\frac{\hat{v}_{\mathfrak{u}}}{\sqrt{n}},{\widetilde{\mathrm{S}}}_{\mathfrak{u},n}^{\textup{closed}} + z_{1- \frac{\alpha}{2}}\frac{\hat{v}_{\mathfrak{u}}}{\sqrt{n}}\right],
$$
with $z_{1- \frac{\alpha}{2}}$ is the $1-\frac{\alpha}{2}$ quantile of the $\mathcal{N}(0,1)$ distribution;
\item[(ii)] for any consistent estimator $\hat{v}_{\sim i}$ of $v_{\sim i}$, where $v_{\sim i}^2$ is defined by (\ref{limitvar}) with $\sim i = \{1, \ldots i-1, i+1, \ldots , p\}$,
an asymptotic CI of level $1- \alpha$ for $\mathrm{S}_{\{i\}}^{\textup{tot}}$ 
is given by
$$\left[ {\widetilde{\mathrm{S}}}_{\{i\},n}^{\textup{tot}} - z_{1- \frac{\alpha}{2}}\frac{\hat{v}_{\sim i}}{\sqrt{n}},{\widetilde{\mathrm{S}}}_{\{i\},n}^{\textup{tot}} + z_{1- \frac{\alpha}{2}}\frac{\hat{v}_{\sim i}}{\sqrt{n}}\right],
$$
with $z_{1- \frac{\alpha}{2}}$ is the $1-\frac{\alpha}{2}$ quantile of the $\mathcal{N}(0,1)$ distribution.
\end{enumerate}
\end{corollary}
\begin{proofof}{{\em Corollary \ref{metmet}}.} Let $\mathfrak{u} \subset \{1 ,\ldots , p\}$. From Item (1) of Theorem 3.4 in \cite{janon:2013}, we know that
$$\sqrt{n}\left(\widetilde{\mathrm{S}}^{\textup{closed}}_{\mathfrak{u},n}-{\mathrm{S}}^{\textup{closed}}_{\mathfrak{u}}\right) \xrightarrow[n \rightarrow + \infty]{\mathcal{D}}\mathcal{N}\left(0,v_{\mathfrak{u}}^2\right)$$
as $n \, \|\tilde{f}_n-f\|_{\infty} \xrightarrow[n \rightarrow + \infty]{}0$, where $v_{\mathfrak{u}}^2$ is the limit variance in Theorem \ref{ICasymp}. Then applying Slutsky's Lemma, it yields the result in Item \textit{(i)}. 

Item \textit{(ii)} for $\mathrm{S}_{\{i\}}^{\textup{tot}}$ is obtained by noting that $\mathrm{S}_{\{i\}}^{\textup{tot}}=1-\mathrm{S}_{\sim i}^{\textup{closed}}$.
\end{proofof}

}

\subsection{Discretized and reduced model}\label{redpart}
\label{sub:discret}
\subsubsection{Discretization of the model}
We use an implicit upwind scheme to discretize \eqref{linear:SW:charac}.  We denote by $k=1,\ldots,N_{t}$ the discrete time index, by $i=1,\ldots,N_{x}$ the space index (where $N_{x}$ and $N_{t}$ are the numbers of discretization points in space and time respectively), and by $\xi_j^{i,k}$ an approximation of $\xi_j$ ($j=1,2$) at the $i$th point of the uniform space grid on $[0,L]$ with $N_{x}$ points and at the $k$th timestep. We denote by $\Delta x = L/N_{x}$ and $\Delta t = T^\star/N_{t}$ the space and time steps, respectively.

Using a classical upwind scheme, we get the following approximations:
\[  \partial_t \xi_1 + \lambda_1 \partial_x \xi_1 \approx \left( \frac{\xi_1^{i,k+1}-\xi_1^{i,k}}{\Delta t} + \lambda_1 \frac{\xi_1^{i,k+1}-\xi_1^{i-1,k+1}}{\Delta x} \right)_{i,k} \]
and:
\[  \partial_t \xi_2 - \lambda_2 \partial_x \xi_2 \approx \left( \frac{\xi_2^{i,k+1}-\xi_2^{i,k}}{\Delta t} - \lambda_2 \frac{\xi_2^{i+1,k+1}-\xi_2^{i,k+1}}{\Delta x} \right)_{i,k}, \]
which give, when combined to \eqref{linear:SW:charac}, the following implicit recurrence (in $k$) relations:
\[ \left( \frac{1}{\Delta t}+\frac{\lambda_1}{\Delta x}+\gamma \right) \xi_1^{i,k+1} - \frac{\lambda_1}{\Delta x} \xi_1^{i-1,k+1} + \delta \xi_2^{i,k+1} = \frac{\xi_1^{i,k}}{\Delta t} \]
\[ \left( \frac{1}{\Delta t}-\frac{\lambda_2}{\Delta x}+\delta \right) \xi_2^{i,k+1} + \frac{\lambda_2}{\Delta x} \xi_2^{i,k+1} + \gamma \xi_1^{i,k+1} = \frac{\xi_2^{i,k}}{\Delta t}. \]
The boundary conditions \eqref{bound:expmodel1} and \eqref{bound:expmodel2} can readily be incorporated so as to write, at each time step, a linear system of equations that has to be solved so as to find $\xi_1^{k+1}$ and $\xi_2^{k+1}$ from $\xi_1^k$ and $\xi_2^k$. 

\subsubsection{Discrete output}
An approximation of the {\it to-be-controlled output} can then be obtained by discretizing the double integral in \eqref{e:defsortie}:
\[ f_{\text{discrete}}(\mu) = \sqrt{ \sum_{i=1}^{N_{x}} \sum_{k=0}^{N_{t}} \left(\xi_1^{i,k}\right)^2 + \left(\xi_2^{i,k}\right)^2 } \]
Notice that we do not include the $1/(N_{x}\times(N_{t}+1))$ normalization factor, as this will not change the Sobol indices.

\subsubsection{Model reduction and error bound}
As the experimental model will have to be numerically solved for a large number of parameter values, we use the \emph{reduced basis} technique (see e.g., \cite{nguyen2005reduced}, and \cite{janon2012goal} for space-time reduction) so as to accelerate the resolutions of the above mentioned systems. We use a space-time reduced basis approach, which can be summed up as follows: we introduce the  vector $\xi = (\xi_1^0, \xi_2^0, \xi_1^1, \xi_2^1, \ldots, \xi_1^{N_{t}}, \xi_2^{N_{t}})$ in $\R^\mathcal N$ (with $\mathcal N = 2 \times (N_{t}+1) \times N_{x}$), which is the solution of the large dimensional problem
$A(\mu) \xi = b(\mu)$, 
where $A(\mu)$ and $b(\mu)$ are appropriate functions of the true parameter values $\mu$.
One can check that $A(\mu)$ and $b(\mu)$ satisfy the so-called affine decomposition hypothesis (as described in \cite{nguyen2005reduced}), hence the classical reduced basis algorithms (op. cit.) can be readily applied to this system, leading to a reduced solution $\tilde\xi$ of the linear system:
\[ \tilde A(\mu) \tilde\xi = \tilde b(\mu), \]
where $\tilde A(\mu)$ is a square matrix of dimension ${m} \ll \mathcal N$, and $\tilde\xi, \tilde b(\mu) \in \R^{m}$. An approximation of $\xi(\mu)$ can then be recovered from $\tilde \xi(\mu)$ from the reduced basis map $Z$:
\begin{equation}\label{e:rbapprox} \xi(\mu) \approx Z \tilde \xi(\mu), \end{equation}
where $Z$ is a suitable matrix with {$m$} columns and $\mathcal N$ lines (see \cite{janon2012goal}). The $Z$ matrix will be assumed, without loss of generality, to have unit and orthogonal column vectors, which are the vectors of the reduced basis. 

We define the standard Euclidean norm: $\norm{x}=\sqrt{x^T x}$. Since the discrete output can be written as
\begin{equation}
\label{eq:f:discrete}
f_{\text{discrete}}(\mu) = \norm{\xi(\mu)}, \end{equation}
the reduced approximation of the output reads:
\begin{equation}
\label{eq:ftilde:discrete}
\tilde f_{\text{discrete}}(\mu) = \norm{ \tilde\xi(\mu) }.\end{equation}

A fully-computable error bound between the solution of the reduced model and the solution of the model is available (op. cit.), which means that the error in \eqref{e:rbapprox} can be quantified. This error bound is given in the theorem below:
\begin{theorem}\label{theo2}\cite{nguyen2005reduced}
For any $\mu$, we have:
\begin{equation}\label{eq:ineq:theom2}
\norm{\xi(\mu) - Z\tilde\xi(\mu)} \leq \frac{\rho(\mu)}{\alpha(\mu)}, \end{equation}
where
\[ \rho(\mu) = \norm{ A(\mu) Z \tilde\xi(\mu) - b(\mu) } \]
and $\alpha(\mu)$ is any real number which satisfies:
\begin{equation}\label{e:defalpha} 0 < \alpha(\mu) \leq \inf_{v, \norm{v}=1} |v^T A(\mu) v|. \end{equation}
\end{theorem}

The papers \cite{nguyen2005reduced} and \cite{huynh2007successive} give efficient methods to compute $\rho(\mu)$, and an $\alpha(\mu)$ which satisfies \eqref{e:defalpha}, respectively.

\begin{corollary}\label{corbound}
For any $\mu$, we have:
\[ \abs{ f_{\text{discrete}}(\mu) - \tilde f_{\text{discrete}}(\mu) } \leq \frac{\rho(\mu)}{\alpha(\mu)}, \]
where $\rho(\mu)$ and $\alpha(\mu)$ are as in Theorem \ref{theo2}.
\end{corollary}
\begin{proofof}{{\em Corollary \ref{corbound}}}
The proof follows from the reverse triangle inequality, (\ref{eq:f:discrete}), (\ref{eq:ftilde:discrete}) and (\ref{eq:ineq:theom2}) in Theorem \ref{theo2}.
\end{proofof}
{\begin{remark}
The bound in Corollary \ref{corbound} depends on the size $m$ of the reduced basis. Hence, we are in the framework of Section \ref{sobolEst}, where our ``true'' output is $Y=f_{\text{discrete}}(\mu)$, and where the metamodel output is $\tilde Y = \tilde f_{discrete,m}(\mu)$. The upper bound for the metamodel error is given by 
$$\delta_m(\mu)=\abs{ f_{\text{discrete}}(\mu) - \tilde f_{discrete,m}(\mu) } \leq \frac{\rho_m(\mu)}{\alpha(\mu)}\, .$$ In the numerical experiments, we will calibrate the reduced basis size $m$ with respect to the Monte-Carlo sample size $n$, in order that $\delta_m(\mu)$ satisfy the assumptions of Corollary \ref{metmet} as the Monte-Carlo sample size $n$ grows to infinity. Note that the calibration of $m$ does not require the evaluation of the true model $f$, just of the error bound, for which efficient algorithms are available, see e.g. \cite{nguyen2005reduced} or \cite{huynh2007successive}.
\end{remark}}

\section{Numerical results}\label{numeric}

\subsection{Parameters}
For the numerical implementations, we have chosen the channel length $L = 250$m, and the time horizon $T^\star = 75$s. 
The parameters $k_0$ and $k_L$ have been fixed to $0.6$ and $0.7$, satisfying condition (\ref{stab:cond:bastin}) for the stability of the nominal closed-loop system.
The discretization parameters were set to $\Delta t = 5$s, $\Delta x = 5$m, and for the reduction, a reduced basis 
obtained from proper orthogonal decomposition \cite{sirovich1987turbulence}
 of size $m$ (to determine, see the next Subsection), obtained from a snapshot of size $100$.
For the estimation of Sobol' indices, the Monte-Carlo sample size was fixed equal to $n=30000$. We provide asymptotic CI
of level $0.95$. 

Note that the length of a CI is a direct measurement of the estimation precision.

\subsection{Calibration of $m$}

In this section, we see how to choose the reduced basis size $m$ depending on the Monte-Carlo sample size $n$, in order to get asymptotic confidence intervals for the sensitivity indices. The idea is to fix $m$ as a function of $n$ so as to satisfy the condition of Corollary \ref{metmet}. We first need to estimate the variance of the error:
\[ \delta(\mu,m) = \tilde f_m(\mu) - f(\mu) \]
as a function of $m$. To do so, we compute samples of $\delta$'s: \[ \mathcal S_m = \{ \delta(\mu,m), \mu\in\Xi \} \]
where $\Xi$ is a random sample of size 1000, and $m=3, 4, \ldots, 14$, and we estimate $\Var  \delta(\mu,m)$ 
by its empirical estimator $\widehat{\Var\delta(\mu,m)}$. The result is given in Figure \ref{f:Bench_reduced_sortie}.

\begin{figure}[ht]
\centerline{\includegraphics[scale=0.35,angle=-90]{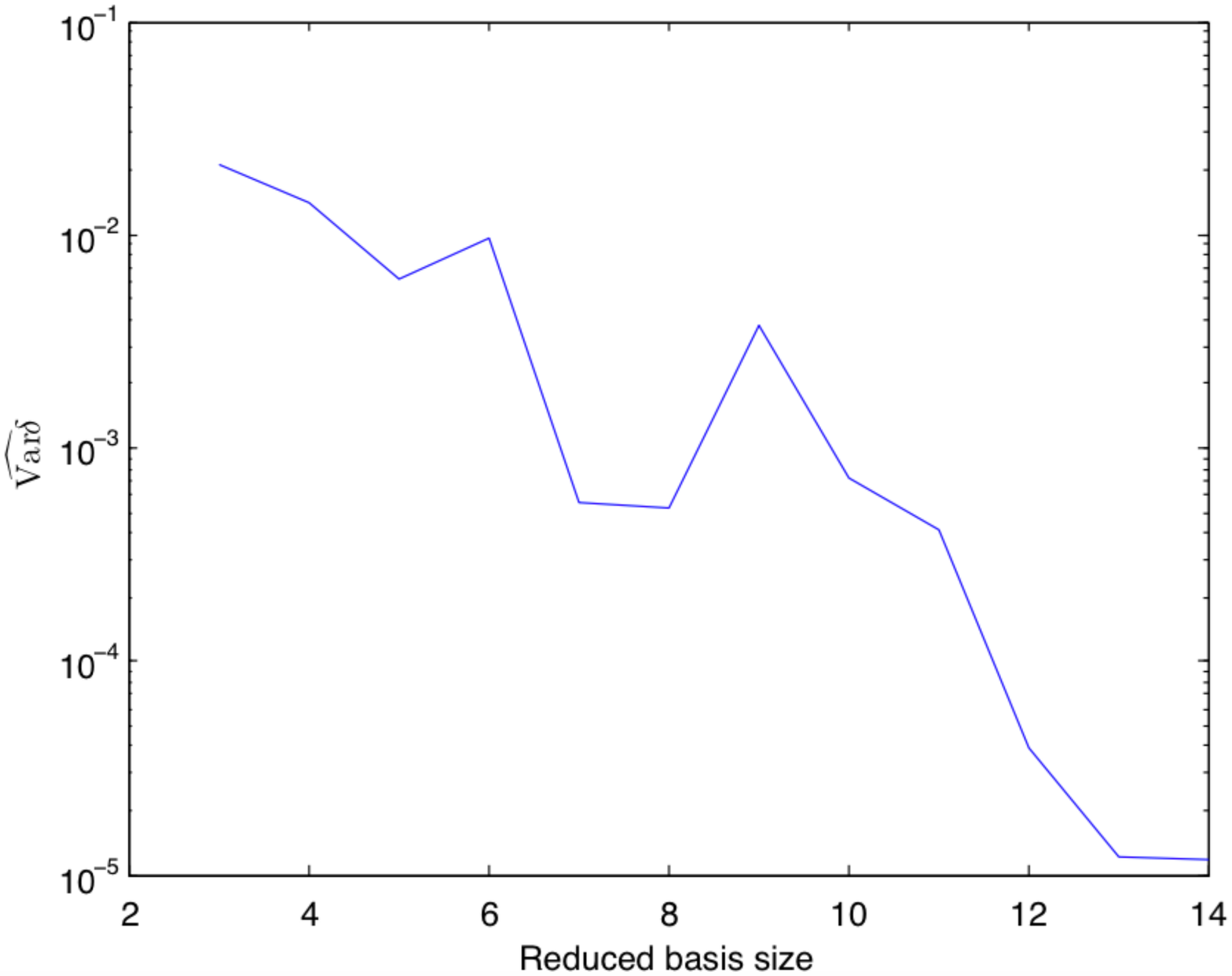}}
\caption{Plot of $\ln \widehat{\Var\delta(\mu,m)}$ as a function of the reduced basis size $m$.}\label{f:Bench_reduced_sortie}
\end{figure}

The plot suggests of a log-linear regression of $\Var\delta(\mu,m)$ 
as a function of $m$. We suppose that the following approximation holds:
\begin{equation}
\label{e:approxregr}
\Var \left( \tilde f_{discrete,m}(\mu) - f_{discrete}(\mu) \right) \approx c q^m
\end{equation}
where $c = 0.2414 $, $q = 0.5070$ (obtained by least-square fitting).

Set $m=m(n)$, and $\tilde f_n=\tilde f_{discrete,m(n)}$. We have, for $n \rightarrow +\infty$:
\[ n \Var \left( \tilde f_n(\mu) - f_{discrete}(\mu) \right) \rightarrow 0 \]
as soon as
\[ n c q^{m(n)} \leq 1/\gamma(n) ,\]
for any $\gamma$ function such that $\lim_{n\rightarrow+\infty} \gamma(n)=+\infty$, that is to say
\[ m(n) \geq - \frac{\log(nc\gamma(n))}{\log q}. \]

To illustrate our asymptotic, we take $n=30000$, and $\gamma(n)=\log\log(n)$.
Hence, we choose  $m$ according to:
\[ m=m(n)=- \frac{\log(nc\log(\log(n)))}{\log q} \approx 14. \]

\subsection{Indices estimations}
%
%
 \begin{figure}[htp]
{\includegraphics[scale=0.33]{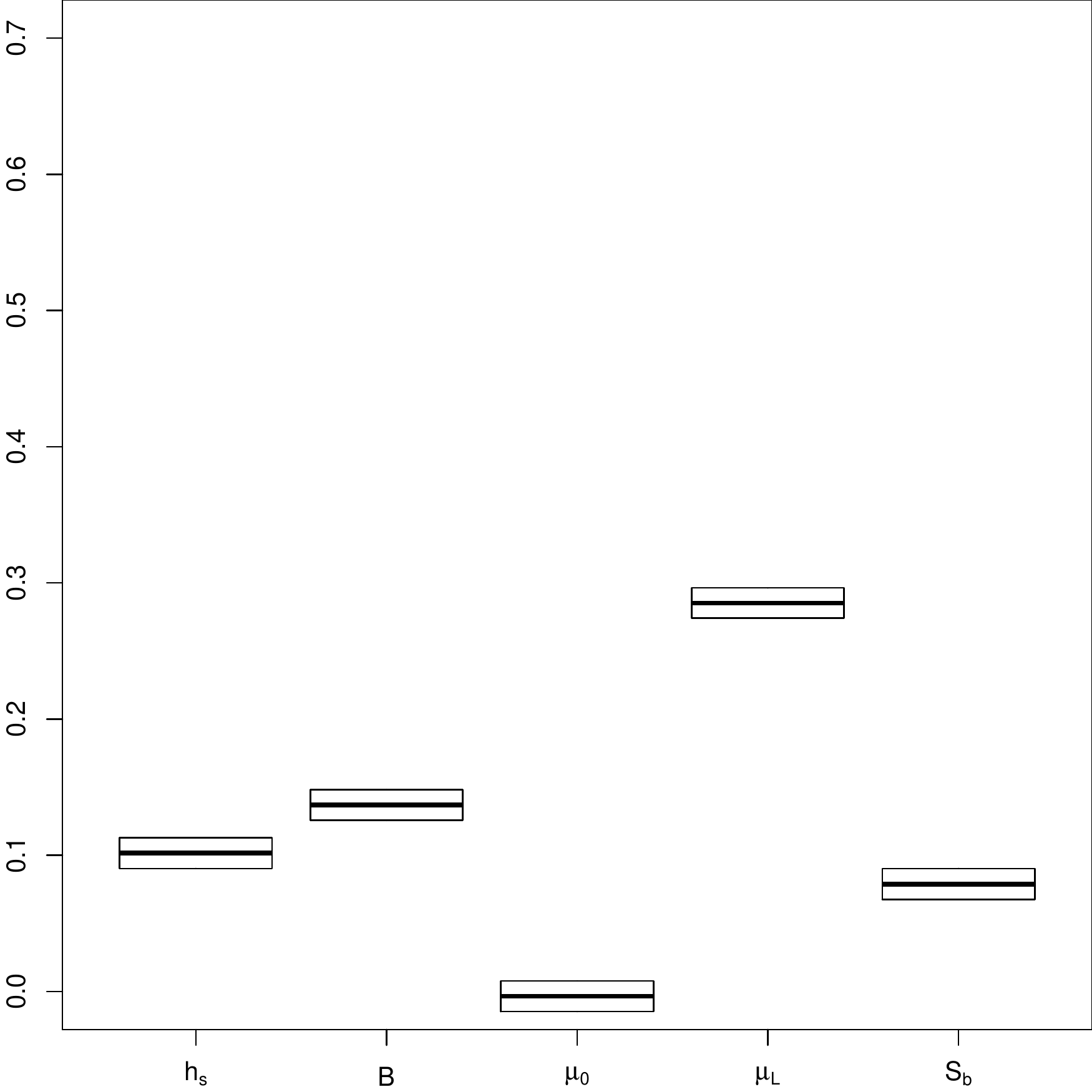}}
{\includegraphics[scale=0.33]{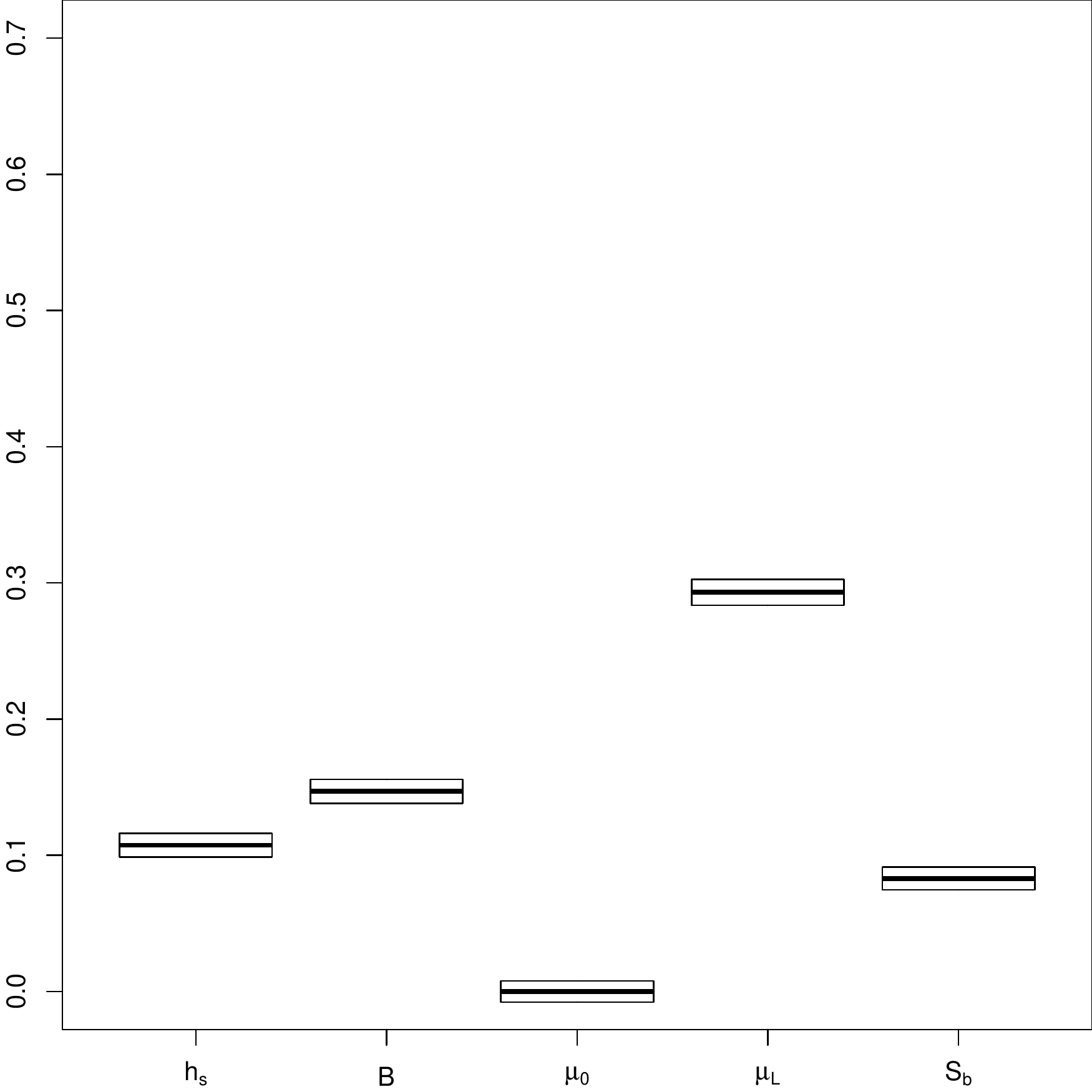}}
\caption{First-order (left) and total (right) Sobol' indices with 95\% CI for for each labelled uncertain parameter.}
\label{fig:graph:order1}
\end{figure}

\begin{figure}[htp]
{\includegraphics[scale=0.33]{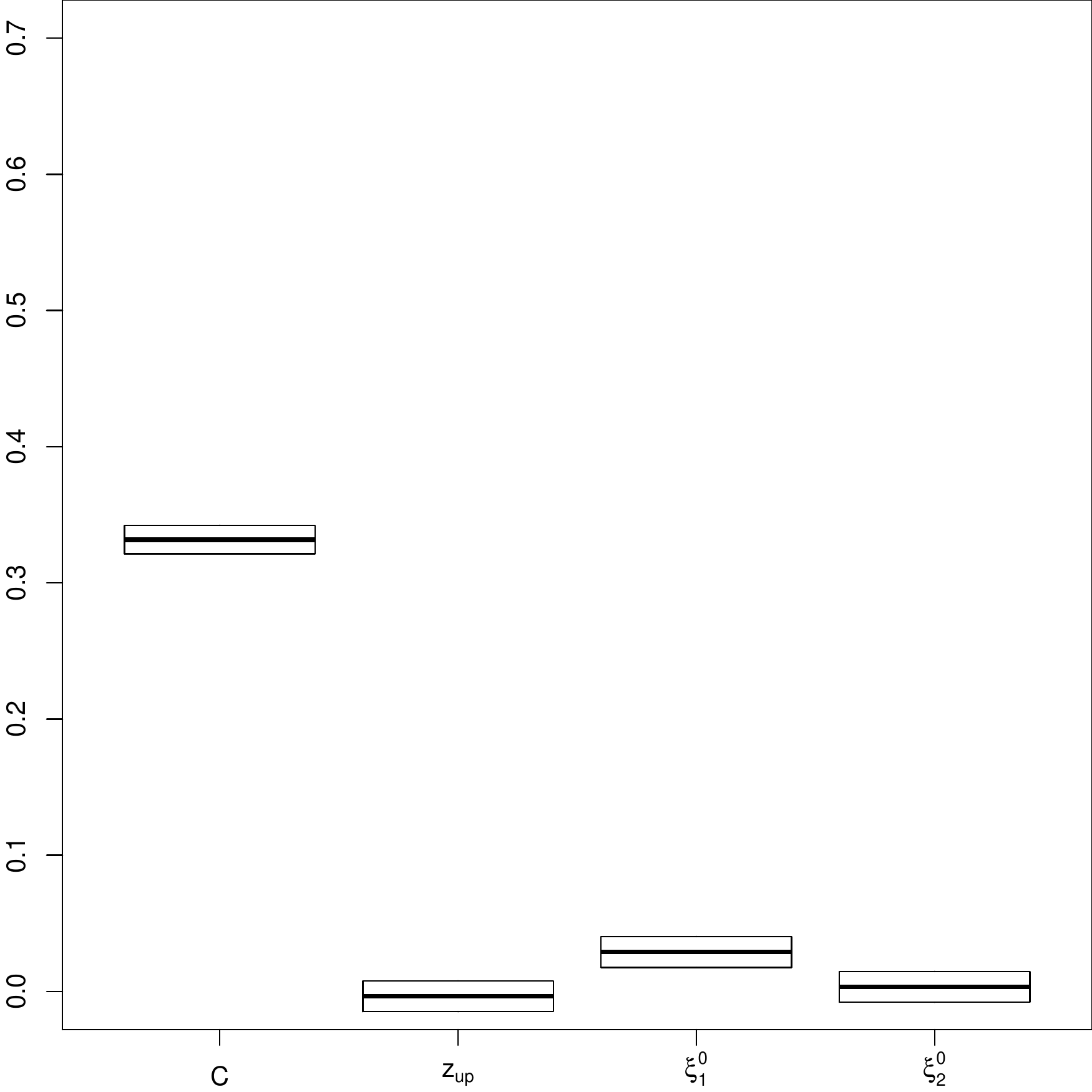}}
{\includegraphics[scale=0.33]{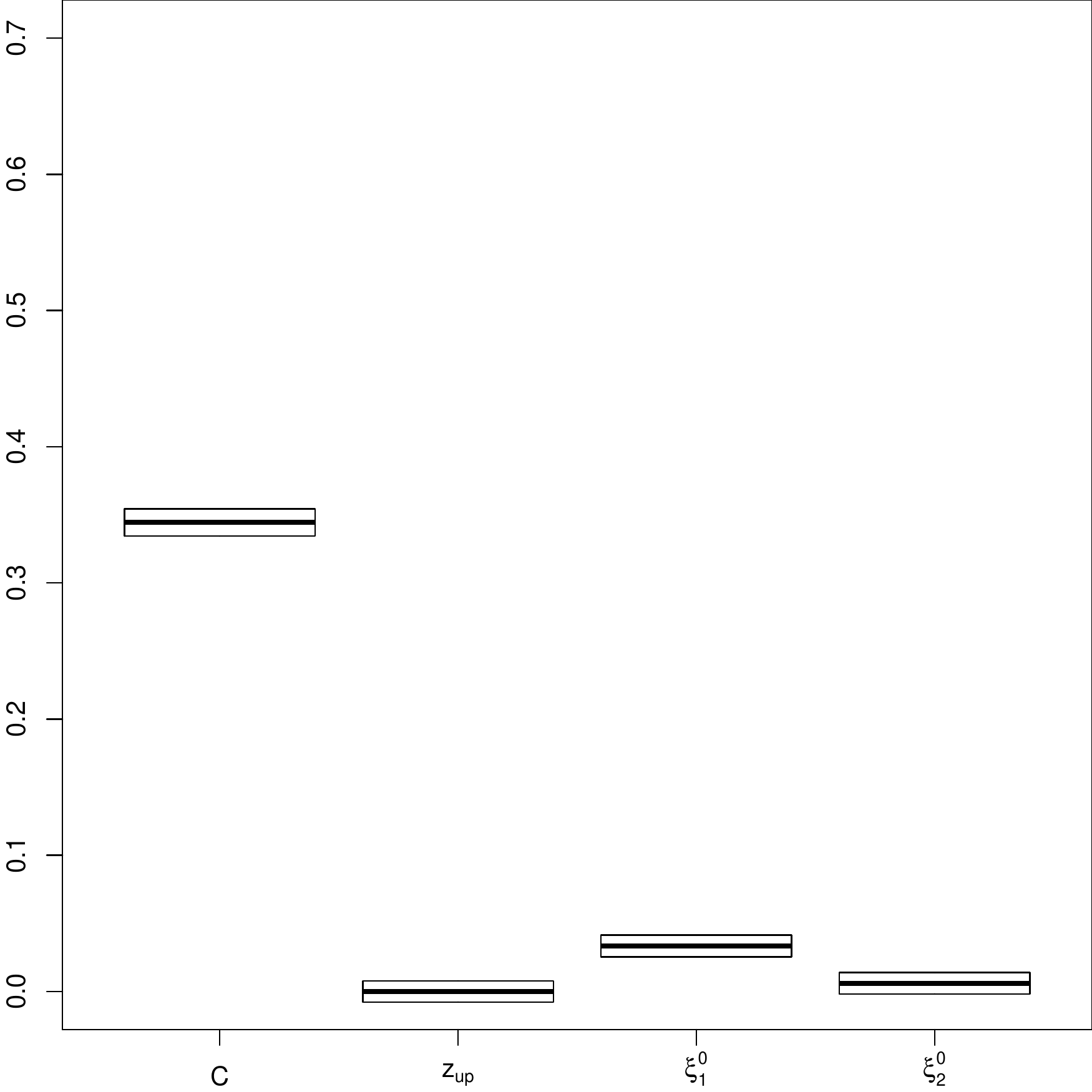}}
\caption{First-order (left) and total (right) Sobol' indices with 95\% CI for  for each labelled uncertain parameter.}
\label{fig:graph:order2}
\end{figure}
On Figures \ref{fig:graph:order1} and \ref{fig:graph:order2} the results read as follows. For each uncertain parameter, the bold centered line on the left part of the figure (respectively on the right part) gives, on the vertical 
axis, the estimation of the first-order (resp. total) Sobol' index. The upper and lower thin lines define the upper and lower bounds of the 95\% CI.
Note that, by definition, the first-order index is smaller than the corresponding total index.
We conclude with a confidence level of $95 \%$ that the influence of the parameters $z_{up}$ and $\mu_0$ are not significant for this  {\em to-be-controlled output}. As expected 
from the choice of  $T^\star$ and the exponential control rate, the initial conditions $(\xi_1^0,\xi_2^0)$ are not significant as well. 
The parameter $C$ is the parameter having the most significant effect, parameters $S_b$ and $B$ are influent as well but to a less extent.

The influence of $\mu_L$ is clearly greater than the influence of $\mu_0$. Therefore as far as the {\em to-be-controlled output} is concerned, a precise knowledge on the value of physical parameter $\mu_L$ is necessary, in contrast to the knowledge of $\mu_0$ which seems to have a smaller impact on the control objective.

For each parameter, as the difference between the total and the first-order indices is not significant, it means that the interactions are negligible. 
Thus the {\em to-be-controlled output} is additive in the parameters.
\subsection{Additive modeling of the output}
The previous results suggested that the to-be-controlled output can be well-approximated by an additive function:
\begin{equation}\label{fAM} f_{discrete} \approx \sum_{x \in \mathcal V} g_x \end{equation}
where $\mathcal V$ is a subset of all the parameters' symbols:
\[ \mathcal V \subset \mathcal V_0=\{ h_s, B, \mu_0, \mu_L, S_b, C, z_{up}, \xi_1^0, \xi_2^0 \}, \]
and the $g_x$ are appropriate univariate functions, that will be fitted using a random iid. sample of 3000 outputs.

There are two possible choices for $\mathcal V$: the ``full'' choice $\mathcal V = \mathcal V_0$, and the ``reduced'' choice:
\[ \mathcal V = \mathcal V_0 \setminus \{ \mu_0, z_{up}, \xi_2^0 \} \]
consisting of all the parameters detected as the most influent by our sensitivity analysis.

The $g_x$ functions can be chosen, either as linear functions (leading to linear regression) or splines (leading to generalized additive modeling (GAM) by splines).

These two-by-two possible choices hence give rise to four possible models. The Akaike information criterion (AIC) is used to compare these four models and select the best one. AIC indeed deals with the trade-off between the goodness of fit of the model and its complexity, see \cite{akaike:1974}. The computations have been made by using the packages \verb<gam< \cite{Rgam} and \verb<stats< \cite{Rstats} packages of the R statistical software \cite{Rstats}. The computed AICs are reported in Table \ref{t:aic}.

\begin{table}[ht]
\begin{center}\begin{tabular}{l|l}
Model & AIC \\ \hline
Full + Linear & -9938 \\ 
Full + Spline &  -10837\\
Reduced + Linear & -9536 \\
Reduced + Spline & -10319
\end{tabular}\end{center}\caption{Fitted Akaike information criterions for the different considered models. The smaller the AIC, the better the model, in terms of bias/variance compromise.}\label{t:aic}
\end{table}

The Table shows the superiority of the GAM (spline) model, compared with the simple linear model. However, using the full set of variables also gives a better AIC than using the reduced one. We hence try:
\begin{equation}\label{e:bonV} \mathcal V = \mathcal V_0 \setminus \{ \mu_0, z_{up} \}, \end{equation}
as $\xi_2^0$ is the ``least negligible'' parameter pointed out by sensitivity analysis. The reported AIC for this choice of $\mathcal V$ is -10847, which is the best AIC obtained. Hence we decide to retain GAM and \eqref{e:bonV} to model $f_{discrete}$.

The  fitted functions $g_x$ are gathered in Figure \ref{f:gams}. We see that all the parameters have a monotonic effect on the output (increasing for $h_s$, $B$, $\mu_L$, decreasing  for the other), and also that all the effects are very close to linear, except for the $B$ parameter, which has a significantly nonlinear (convex) effect.
\begin{figure}
\includegraphics[scale=.35,page=1]{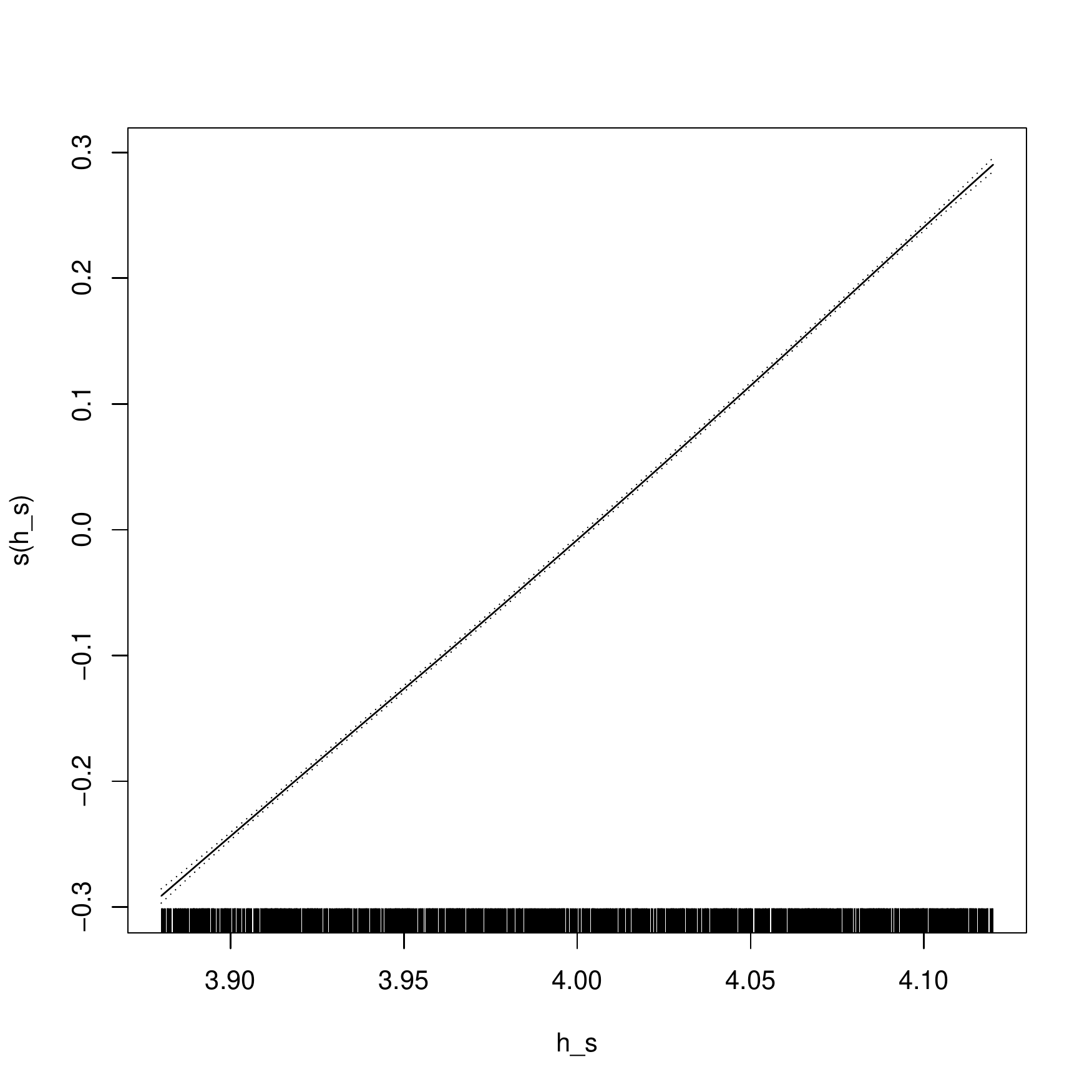}
\includegraphics[scale=.35,page=2]{GAMfits.pdf}
\includegraphics[scale=.35,page=3]{GAMfits.pdf}
\includegraphics[scale=.35,page=4]{GAMfits.pdf}
\includegraphics[scale=.35,page=5]{GAMfits.pdf}
\includegraphics[scale=.35,page=6]{GAMfits.pdf}
\includegraphics[scale=.35,page=7]{GAMfits.pdf}
\label{f:gams}
\caption{The fitted $g_x$ spline functions in \eqref{fAM}, for $x\in\mathcal V$ defined by \eqref{e:bonV}. The dotted lines are the 95\% confidence bands. }
\end{figure}

\section{Conclusion}\label{sec:conc}

In this paper, the global sensitivity analysis has been performed when considering probabilistic distribution functions standing for the uncertainty in some unknown physical parameters. It allows us to describe the impact of the parameters on a {\it to-be-controlled output} in a boundary control problem. 
This boundary control is motivated by an application for the flow control in an open channel where water height and velocity are described by the Shallow Water equations, in presence of friction and bottom slopes. The {\em to-be-controlled output} has been defined as the norm of the state at a a given large time. We deduce from the numerical studies that some parameters have a nonlinear effect (such as $B$) whereas other parameters have a linear effect, such as the parameter $C$ appearing in the friction slope $S_f$, which is usually not well modeled in the literature.
Concerning the physical parameters $\mu_0$ and $\mu_L$, it appears that a precise knowledge of $\mu_L$ is more important than the knowledge of $\mu_0$ (since the latter parameter is less influent on the {\em to-be-controlled output}).

This work lets many research lines open. In particular, it could be interesting to optimize the control parameters  $k_0$ and $k_L$ by adding them
in the global sensitivity analysis. Another open question is to provide the {\em a posteriori}  error estimator for the reduced basis approach with a nonlinear output. This latter study  is currently under  investigation.

\appendix
\section{Technical proofs}\label{sec:calc}

\subsection{Proof of Proposition \ref{def:controls:lemma}}
\label{appendix:A}

\begin{proofof}{{\em Proposition \ref{def:controls:lemma}}.}
Let us first note that it follows from (\ref{tau3}) and $Q=BHV$ that (by omitting the time variable)
$$
B H(0) V(0) = U _0 B \mu _0\sqrt{2g(z_{up}-H(0))}
$$
which is equivalent to 
\begin{equation}
\label{tau3.1}
V(0) = U _0 \frac{1}{H(0)} \mu _0\sqrt{2g(z_{up}-H(0))}\ .
\end{equation}

On the other side the first line of (\ref{bound:linear:SW}) yields with (\ref{charac:coord}) and the definitions of $v$ and $h$, the following
\begin{equation*}\begin{array}{c}
V(0) - V ^\star + ( H (0) - H ^\star ) \sqrt{\frac{g}{H^\star}} \\
= k_0 ( V(0)- V ^\star) - k_0 (H(0) - H ^\star) \sqrt{\frac{g}{H^\star}} 
\ ,\end{array}\end{equation*}
which may be rewritten as
\begin{equation}
\label{bound:linear:SW.1}
V(0) =V ^\star -\frac{1+ k_0}{1- k_0} ( H (0) - H ^\star ) \sqrt{\frac{g}{H^\star}} \ .
\end{equation}
Therefore, with (\ref{tau3.1}) and (\ref{bound:linear:SW.1}), (\ref{tau3}) and the first line of (\ref{bound:linear:SW}) are equivalent as soon as the first control is defined by (\ref{def:U0}).

Let us compute the control $U_L$ in a similar way.
To do that, we first deduce from (\ref{under}) the following
\begin{equation*}
(H(L)- h _ s  - U _L)^3 = \frac{ H (L) ^2V (L) ^2}{2 g  \mu_L ^2}
\end{equation*}
and thus
\begin{equation}
\label{under.2}
V(L)=\frac{\sqrt{2 g \mu_L ^2(H(L)- h _ s  - U _L)^3}}{H (L)}
\end{equation}
Moreover from the second line of (\ref{bound:linear:SW}), with (\ref{charac:coord}) and the expressions of $v$ and $h$, it holds
\begin{equation*}
\begin{array}{c}
V(L) -V ^\star  - ( H( L) - H ^\star) \ds\sqrt{\frac{g }{H^\star}} \\= k_L  (V(L) -V ^\star) + k_L (H(L) -H ^\star)   \ds\sqrt{\frac{g }{H^\star}}\end{array}
\end{equation*}
and also
\begin{equation}
\label{bound:linear:SW.2}
V(L) = V ^\star  + \frac{1+ k_L}{1- k_l} ( H( L) - H ^\star) \sqrt{\frac{g }{H^\star}}
\end{equation}
From (\ref{under.2}) and (\ref{bound:linear:SW.2}) we get that (\ref{under}) and the second line (\ref{bound:linear:SW}) are equivalent as soon as the control $U_L$ is defined by
\begin{equation*}\begin{array}{c}
\sqrt{2 g \mu_L ^2  ( H(L) - h_s - U_L ) ^3} \\
= H (L)(V ^\star  + \frac{1+ k_L}{1- k_L} ( H( L) - H ^\star) \sqrt{\frac{g }{H^\star}}
)
\end{array}\end{equation*}
which is equivalent to (\ref{def:UL}). This concludes the proof of Proposition \ref{def:controls:lemma}.
\end{proofof}

\subsection{Proof of Proposition \ref{def:controlsreal_life:lemma}}
\label{appendix:B}
%
%
\begin{proofof}{{\em Proposition \ref{def:controlsreal_life:lemma}}.}
Let us define
$ \alpha=\frac{1+k_0}{1-k_0}, \;\; \beta=\sqrt{\frac{g}{\hetn}} .$

At $x=0$, we take, from the proof of Proposition \ref{def:controls:lemma}:
\[ H(0) V(0) = U_0 \mu_0 \sqrt{2g(z_{up}-H(0))} \]
with
\[ U_0 = \frac{H(0)(\vetn-\alpha(H(0)-\hetn)\beta)}{\mu_0\sqrt{2g(z_{up,\nom}-H(0))}} \]
that is:
\[ v(0) = \sqrt{\frac{z_{up}-H(0)}{z_{up,\nom}-H(0)}} \left( \vetn - \alpha ( H(0)-\hetn ) \beta \right) - \vet. \]
Now, given the change of variables between $(V,H)$, $(v,h)$,  this equation could be rewritten as a nonlinear relation between $v(0)$ and $h(0)$, hence, as a nonlinear relation between $\xi_1(0)$ and $\xi_2(0)$. However, to keep the resolution simple, we have to linearize this equation. This linearization is made in accordance with the linearization of the Shallow-Water equation for $(H,V)$ near $(H^\star, V^\star)$, hence for $(h,v)$ near the origin. We do the same here; hence by Taylor expansion around $h=0$, we get
$ v(0)= \calA + \calB h(0) + o (h(0)) $
with
\[ \calA = \frac{\mu_0}{\mu_{0,\nom}}\sqrt{ \frac{\het-z_{up} }{\het-z_{up,\nom}}}\left( \vetn-\alpha\beta(\het-\hetn) \right) - \vet, \]
and 
$$
\begin{array}{c}\calB=\\\frac{\mu_0}{\mu_{0,\nom}}\sqrt{\frac{\het-z_{up} }{\het-z_{up,\nom}}}\left( - \alpha\beta + \frac{\left(z_{up}-z_{up,\nom}\right)\left( \vetn - \alpha \beta (\het-\hetn) \right) }{2(\het-z_{up})(\het-z_{up,\nom})}\right). \end{array}$$

Similarly, at $x=L$, we have
\[ V(L) = \frac{\sqrt{2g } \mu_L}{H(L)} \left( e_h + \left( \frac{H(L)(\vetn+\alpha_L \beta e_H)}{\sqrt{2g} \mu_{L}} \right)^{2/3} \right)^{3/2} \]
where $\alpha_L=\frac{1+k_L}{1-k_L}$, and
\[ e_h = h_{s,\nom} - h_s, \;\; e_H = H(L)- H^\star_\nom. \]
Therefore,
\[ v(L) = \calC + \calD h(L) + o(h(L)), \]
with
\[ \begin{array}{rcl}\calC &=&{\frac{\sqrt{2g}}{4}}\mu_L \left(2 h_{s,\nom}-2h_s\right.
\\&&\left.+2^{\frac{2}{3}} \left(-H^\star \frac{\left(-V^\star_{\nom}+\alpha_L \beta_{\nom} \times H^\star_{\nom}\right)}{\left(\sqrt{g}\mu_{L,\nom}\right)}\right)^{\frac{2}{3}}\right) 
\\
&&\frac{\sqrt{\left(4 h_{s,\nom}-4 h_s+2\times 2^{\frac{2}{3}} \left(-H^\star\times \frac{\left(-V^\star_{\nom}+\alpha_L \beta_{\nom}  H^\star_{\nom}\right)}{\left(\sqrt{g} \mu_{L,\nom}\right)}\right)^{\frac{2}{3}}\right)}}{H^\star}\end{array}\]

\[ \begin{array}{rcl}\calD &=& -{\frac{1}{2}} \left(2^{\frac{1}{6}} \left(-H^\star \frac{\left(-V^\star_{\nom}+\alpha_L  \beta_{\nom} H^\star_{\nom}\right)}{\left(\sqrt{g} \mu_{L,\nom}\right)}\right)^{\frac{2}{3}} \alpha_L  \beta_{\nom} H^\star
\right.\\&&
\left.-\sqrt{2} V^\star_{\nom} h_{s,\nom}+\sqrt{2} V^\star_{\nom} h_s+\sqrt{2} \alpha_L  \beta_{\nom}  H^\star_{\nom} h_{s,\nom}
\right.\\
&&\left.-\sqrt{2}\alpha_L  \beta_{\nom}  H^\star_{\nom} h_s\right)
\\
&&\sqrt{\left(4 h_{s,\nom}-4h_s+2 \times 2^{\frac{2}{3}} \left(-H^\star\frac{\left(-V^\star_{\nom}+\alpha_L  \beta_{\nom} \times H^\star_{\nom}\right)}{\left(\sqrt{g}\mu_{L,\nom}\right)}\right)^{\frac{2}{3}}\right)}
\\&& \mu_L  \frac{\sqrt{g}}{\left(\left(-V^\star_{\nom}+\alpha_L \beta_{\nom} H^\star_{\nom}\right)H^{\star 2} \right)}\end{array}
 \]
Thus, the linearized boundary relations satisfied by $\xi_1$ and $\xi_2$ in the real-life model are given by (\ref{bound:expmodel1}) and (\ref{bound:expmodel2}).
\end{proofof}

\bibliographystyle{plain}      
\bibliography{JNPP-mar14}
\end{document}